\documentclass{amsart}

\usepackage{amsmath}
\usepackage{amsfonts}
\usepackage{amssymb,enumerate}
\usepackage{amsthm}
\usepackage[all]{xy}

\newtheorem{lem}{Lemma}[section]
\newtheorem{cor}[lem]{Corollary}
\newtheorem{prop}[lem]{Proposition}

\newtheorem{Defn}[lem]{Definition}
\newtheorem{Ex}[lem]{Example}
\newtheorem{Question}[lem]{Question}
\newtheorem{Property}[lem]{Property}
\newtheorem{Properties}[lem]{Properties}
\newtheorem{Discussion}[lem]{Remark}
\newtheorem{Construction}[lem]{Construction}
\newtheorem{Subprops}{}[lem]
\newtheorem{Para}[lem]{}

\newtheorem{intthm}{Theorem}

\newenvironment{ex}{\begin{Ex}\rm}{\end{Ex}}

\newenvironment{para}{\begin{Para}\rm}{\end{Para}}
\newenvironment{disc}{\begin{Discussion}\rm}{\end{Discussion}}

\newcommand{\comp}[1]{\widehat{#1}}
\newcommand{\Comp}[2]{\widehat{#1}^{\ideal{#2}}}
\newcommand{\ideal}[1]{\mathfrak{#1}}
\newcommand{\m}{\ideal{m}}
\newcommand{\fa}{\ideal{a}}
\newcommand{\fb}{\ideal{b}}
\newcommand{\n}{\ideal{n}}

\newcommand{\fp}{\ideal{p}}

\newcommand{\pd}{\operatorname{pd}}

\newcommand{\ext}{\operatorname{Ext}}

\newcommand{\rhom}{\mathbf{R}\mathrm{Hom}}      
\newcommand{\lotimes}{\otimes^{\mathbf{L}}}

\newcommand{\HH}{\operatorname{H}}
\newcommand{\Hom}{\operatorname{Hom}}

\newcommand{\vf}{\varphi}

\newcommand{\D}{\mathsf{D}}

\newcommand{\ann}{\operatorname{Ann}}

\newcommand{\xra}{\xrightarrow}

\newcommand{\shift}{\mathsf{\Sigma}}

\newcommand{\M}[2]{\mathbf{C}^{\ideal{#1}}(#2)}
\newcommand{\image}{\operatorname{Im}}

\newcommand{\fc}[2]{\operatorname{C}^{\ideal{#1}}_{#2}}
\newcommand{\LL}[2]{\mathbf{L}\Lambda^{\ideal{#1}}(#2)}
\newcommand{\RG}[2]{\mathbf{R}\Gamma_{\ideal{#1}}(#2)}

\newcommand{\colim}{\operatorname{colim}}
\newcommand{\llim}{\operatorname{lim}}
\newcommand{\ass}{\operatorname{Ass}}
\newcommand{\supp}{\operatorname{Supp}}

\newcommand{\cpln}[2]{\varepsilon^{\ideal{#1}}_{#2}}
\newcommand{\inc}[2]{i^{\ideal{#1}}_{#2}}
\newcommand{\eff}[2]{f^{\ideal{#1}}_{#2}}
\newcommand{\eval}[2]{g^{\ideal{#1}}_{#2}}
\newcommand{\deval}[2]{h^{\ideal{#1}}_{#2}}
\newcommand{\ceval}[2]{k^{\ideal{#1}}_{#2}}
\newcommand{\adj}{\theta^{}}

\renewcommand{\leq}{\leqslant}
\renewcommand{\geq}{\geqslant}

\begin{document}

\bibliographystyle{amsplain}

\author{Anders J.~Frankild}
\address{University of Copenhagen, Institute for Mathematical 
Sciences, Department of Mathematics, 
Universitetsparken 5, 2100 K\o benhavn, Denmark}
\email{frankild@math.ku.dk}
\urladdr{http://www.math.ku.dk/~frankild/}

\author{Sean Sather-Wagstaff}
\address{Department of Mathematics, California State University, Dominguez Hills, 
1000 E.~Victoria St., Carson, CA 90747 USA}
\curraddr{Department of Mathematical Sciences, Kent State University,
  Kent OH
  44242}
\email{sather@math.kent.edu}
\urladdr{http://www.math.kent.edu/~sather}

\thanks{This research 
was conducted while 
AJF had a Steno Stipend from the Danish Research Council.}

\title{Detecting completeness from ext-vanishing}

\dedicatory{Dedicated to Lex Remington}

\begin{abstract}
Motivated by work of C.~U.~Jensen, R.-O.~Buchweitz, and H.~Flenner, we prove the following result.
Let $R$ be a commutative noetherian ring and $\fa$
an ideal
in the Jacobson radical of $R$.  
Let $\Comp{R}{a}$ be the $\fa$-adic completion of $R$.
If $M$ is a finitely generated $R$-module 
such that $\ext^i_R(\Comp{R}{a},M)=0$ for all $i\neq 0$,
then $M$
is $\fa$-adically complete.
\end{abstract}
\maketitle

\section*{Introduction} 

A result of Jensen~\cite[(8.1)]{jensen:lfdlatm} characterizes the completeness property of a semilocal
ring in terms of Ext-vanishing:  If $R$ is a commutative noetherian ring, 
then it is a finite product of complete local rings if and only if 
$\ext^i_R(B,M)=0$ for $i\neq 0$ whenever $B$ is flat and $M$ is finitely generated over $R$.
In their investigation of Hochschild homology,
Buchweitz and Flenner~\cite[(2.3)]{buchweitz:psrp} recover one implication of
the local case of this result:
Let $R$ be a ring and $\m\subset R$ a maximal ideal; 
if $M$ is an $\m$-adically complete $R$-module,
then $\ext^i_R(B,M)=0$ for all $i\neq 0$ and each flat $R$-module $B$;
see also~\cite[(3.7)]{frankild:volh} for the local case.

In this paper, we investigate  converses to the Buchweitz-Flenner 
result:  If $M$ is an $R$-module such that
$\ext^i_R(B,M)=0$ for all $i\neq 0$ and each flat $R$-module $B$, must $M$ be $\m$-adically
complete?  One  sees readily that
this need not be the case when $M$ is not finitely generated. If $R$ is a 
local domain with $\dim(R)>0$ and $M$ is the quotient field of $R$, then 
$M$ is not $\m$-adically complete.  However, $M$
is injective so $\ext^i_R(B,M)=0$ for all $i\neq 0$ and each $R$-module $B$.

The following result is proved in~\ref{proofA}.
When $M$ finitely generated,
it shows that the completeness of $M$ can be 
ascertained from the vanishing of the Ext-modules against a single flat module, 
namely $\comp{R}$.

\begin{intthm} \label{thma}
Let $R$ be a commutative noetherian ring and $\fa$
an ideal
in the Jacobson radical of $R$.  
Let $\Comp{R}{a}$ be the $\fa$-adic completion of $R$
and let $M$ be a finitely generated $R$-module.
The following conditions are equivalent.
\begin{enumerate}[\quad\rm(i)]
\item \label{thmaitem1} 
$M$ is $\fa$-adically complete.
\item \label{thmaitem2} 
$\ext^i_R(\Comp{R}{a},M)=0$ for all $i\neq 0$.
\item \label{thmaitem3} 
$\ext^i_R(\Comp{R}{a},M)=0$ for all $i =1,\ldots,\dim_R(M)$.
\end{enumerate}
\end{intthm}

As a consequence of this theorem we obtain the following two results.
The first is proved in~\ref{proofreferee2}, and the second is contained in
Corollary~\ref{cor03}.

\begin{intthm} \label{referee2}
The ring $R$ is $\fa$-adically complete
if and only if the completion
$\Comp{R}{a}$ is module-finite over $R$.
\end{intthm}

\begin{intthm} \label{thmb}
Let $M,N$ be  finitely generated $R$-modules
and $t$  an integer such that $\ext_R^i(N,M)=0$ for each $i<t$.
If $\ext^i_R(\Comp{N}{a},M)=0$ for each $i\neq t$,
then $\ext^i_R(N,M)=0$ for each $i\neq t$ and 
$\ext_R^{t}(N,M)$ is $\fa$-adically complete.
\end{intthm}

To prove these results, we employ a combination of classical module-theory and
derived category techniques.  Preliminary module-theoretic results are presented in
Section~\ref{sec1}.  Requisite derived category notions are discussed in Section~\ref{sec0}.

\section{Analytic conductor submodules} \label{sec1}

\noindent\emph{Throughout
this work, $R$ is a commutative noetherian ring and $\fa$ is an ideal
contained in the Jacobson radical of $R$.}

\begin{lem} \label{lem01}
If  $M$ is a finitely generated $R$-module,
then $M$ admits a unique maximal $\fa$-adically complete submodule $\fc{a}{M}$.
\end{lem}

\begin{proof}
Let $\M{a}{M}$ denote the collection of $\fa$-adically complete submodules of $M$
which is nonempty because it contains the zero submodule. Since $M$ is noetherian,
this collection contains maximal elements,
each of which is finitely generated.  Let $N,N'\in\M{a}{M}$ be maximal elements and 
suppose that $N\neq N'$.  By maximality, one has $N\not\subseteq N'$ and so
$N\subsetneq N+N'$.  In particular, $N+N'$ is not $\fa$-adically complete.
However, the module $N\oplus N'$ is finitely generated and $\fa$-adically complete.
Hence, the homomorphic image $N+N'$
of $N\oplus N'$ is $\fa$-adically complete, a contradiction.
Thus, $N=N'$ and the maximal element of 
$\M{a}{M}$ is unique.
\end{proof}

The submodule $\fc{a}{M}$ is the \emph{analytic conductor} of $M$ with respect to $\fa$. 
It is the largest $R$-submodule of $M$ that is also an $\Comp{R}{a}$-module.
Before presenting an important property of $\fc{a}{M}$ for this work,
we introduce some frequently used maps.

\begin{para} \label{maps}
Let $M$ be an $R$-module.
The map $\eval{a}{M}\colon \Hom_R(\Comp{R}{a},M)\to M$ is given by $\eval{a}{M}(\vf)=\vf(1)$,
and $\cpln{a}{M}\colon M\to \Comp{M}{a}$ is the natural inclusion.
Assume now that $M$ is 
finitely generated, so that $\fc{a}{M}$ is defined.
Let $\inc{a}{M}\colon \fc{a}{M}\to M$
denote the natural inclusion.  
The map $\eff{a}{M}\colon \fc{a}{M}\to\Hom_R(\Comp{R}{a},M)$
is given by $\eff{a}{M}(m)(r)=rm$.  

The next result yields a well-defined  map
$\ceval{a}{M}\colon \Hom_R(\Comp{R}{a},M)\to \fc{a}{M}$, given by $\ceval{a}{M}(\vf)=\vf(1)$,
such that $\eval{a}{M}=\inc{a}{M}\ceval{a}{M}$.
\end{para}

\begin{lem} \label{lem02}
If  $M$ is a finitely generated $R$-module,
then the natural inclusion 
$\Hom_R(\Comp{R}{a},\inc{a}{M})\colon\Hom_R(\Comp{R}{a},\fc{a}{M})\to\Hom_R(\Comp{R}{a},M)$
is bijective.
\end{lem}

\begin{proof}
By left-exactness of $\Hom_R(\Comp{R}{a},-)$ the given map is injective.
To see that this map is surjective,
fix $\vf\in\Hom_R(\Comp{R}{a},M)$; it suffices to show $\image(\vf)\subseteq\fc{a}{M}$.
The image $\image(\vf)$ is finitely generated over $R$ and a homomorphic image of the
$\fa$-adically complete $R$-module $\Comp{R}{a}$.
Hence,  $\image(\vf)$ is $\fa$-adically complete, and
the desired conclusion follows from
Lemma~\ref{lem01}. 
\end{proof}

\section{Derived local homology and cohomology} \label{sec0}

We work in the derived category $\D(R)$ of 
complexes of $R$-modules, indexed homologically.
References on the subject include~\cite{gelfand:moha, hartshorne:rad}.
A complex $X$ is \emph{homologically bounded to the right} if $\HH_i(X)=0$ for all $i\ll 0$;  
it is \emph{homologically degreewise finite} if $\HH_i(X)$ is finitely generated for each $i$;  
it is \emph{homologically finite} if $\oplus_i\HH_i(X)$  is finitely generated;  and
it is \emph{homologically concentrated in degree $s$} if $\HH_i(X)=0$ for all $i\neq s$.
Isomorphisms in $\D(R)$ are identified by the symbol $\simeq$,
as are quasiisomorphisms in the category  of complexes.
For $X,Y\in\D(R)$ set
$\inf(X)$ and $\sup(X)$ to be the infimum and supremum, respectively, of the set
$\{n\in\mathbf{Z}\mid\HH_n(X)\neq 0\}$.
Let $X\lotimes_R Y$ and $\rhom_R(X,Y)$ denote the left-derived
tensor product and right-derived homomorphism complexes, respectively.

The left-derived local homology and 
right-derived local cohomology functors with support in 
an ideal $\fa$ are denote
$\LL{a}{-}$ and $\RG{a}{-}$, respectively;
see~\cite{lipman:lhcs,greenlees:dfclh}. 
These are computed as follows.  If
$P\xra{\simeq} X\xra{\simeq} J$ are K-projective and K-injective resolutions, 
respectively, as in~\cite{avramov:hdouc,spaltenstein:ruc},
then 
\begin{align*}
\Lambda^{\fa}(-)&=\llim_{n}(R/\fa^n\otimes_R -) &
\Gamma_{\fa}(-)&=\colim_{n}\Hom_R(R/\fa^n,-)\\
\LL{a}{X}&= \Lambda^{\fa}(P) &
\RG{a}{X}&=\Gamma_{\fa}(J).
\end{align*}
Note that the functor $\Gamma_{\fa}(-)$ is left-exact
while $\Lambda^{\fa}(-)$ is neither left- nor right-exact.

\begin{para} \label{props01}
Here is a catalog of properties of $\LL{a}{-}$ and $\RG{a}{-}$ that we will utilize.
\begin{enumerate}[\rm(a)]
\item\label{item01a}
There are natural transformations  of functors on $\D(R)$; \cite[(0.3)$^*$]{lipman:lhcs}.
$$\RG{a}{-}\xra{\gamma} 1_{\D(R)}(-)\xra{\nu}\LL{a}{-}$$
\item 
The following 
are 
equivalences of functors on $\D(R)$;
\cite[Cor.~to (0.3)$^*$]{lipman:lhcs}. \label{item01c}
$$
\LL{a}{\RG{a}{-}}\xra{\LL{a}{\gamma}} \LL{a}{-}
\qquad\text{and}\qquad
\RG{a}{-}\xra{\RG{a}{\nu}} \RG{a}{\LL{a}{-}}
$$
\item\label{item01d}  One has natural equivalences of functors on $\D(R)$; 
\cite[(0.3)]{lipman:lhcs} and~\cite[(3.1.2)]{lipman:llcd}.
$$
\LL{a}{-}\simeq\rhom_R(\RG{a}{R},-)
\qquad\text{and}\qquad
\RG{a}{-}\simeq\RG{a}{R}\lotimes_R -
$$
\item\label{item01e}
(Adjointness) There is a natural equivalence of bifunctors on $\D(R)$
$$\rhom_R(\RG{a}{-},-)\xra[\simeq]{\adj}\rhom_R(-,\LL{a}{-})$$
such that, for all complexes $X$ and $Y$ the next diagram commutes;
\cite[(0.3)]{lipman:lhcs}.
$$\xymatrix{
\rhom_R(X,Y) \ar[d]_{\rhom_R(\gamma_X^{},Y)} \ar[rrd]^{\qquad\rhom_R(-,\nu^{}_Y)} \\
\rhom_R(\RG{a}{X},Y)\ar[rr]^{\adj_{XY}}_{\simeq} &&\rhom_R(X,\LL{a}{Y})
}
$$
In particular, the morphism $\rhom_R(\gamma^{}_X,Y)$
is an isomorphism in $\D(R)$ if and only if
$\rhom_R(-,\nu^{}_Y)$ is so.
\item\label{item01f}
One has a natural equivalence
of functors on 
the full subcategory of $\D(R)$ 
of homologically degreewise finite and bounded to the right complexes;
\cite[(2.8)]{frankild:volh}.
$$\LL{a}{-}\simeq-\lotimes_R\Comp{R}{a}$$
\item\label{item01g} Parts~\eqref{item01c}--\eqref{item01d} 
yield
equivalences 
of (bi)functors on $\D(R)$; see e.g.~\cite[(A.4.22)]{christensen:gd}. 
\begin{gather*}
\RG{a}{R}\lotimes_R\LL{a}{-}\simeq\RG{a}{-}\\
\LL{a}{\rhom_R(-,-)}\simeq\rhom_R(-,\LL{a}{-})
\end{gather*}
\item\label{item01j}
If $X$ is homologically bounded to the right, then
it admits a K-projective resolution $P\xra{\simeq}X$ such that
$X_i=0$ for each $i\leq\inf(X)$, and so
$$\inf(\LL{a}{X})=\inf(\Lambda^{\fa}(P))\geq\inf(P)=\inf(X).$$
\end{enumerate}
\end{para}

We now verify facts about $\LL{a}{-}$ and $\RG{a}{-}$ for the sequel.
Fix $M\in\D(R)$ with K-injective resolution $M\xra{\simeq} J$.  
The map 
$\eval{a}{J}\colon\Hom_R(\Comp{R}{a},J)\to J$
given by $\vf\mapsto\vf(1)$ 
describes a well-defined morphism
$\deval{a}{M}\colon\rhom_R(\Comp{R}{a},M)\to M$
in $\D(R)$.

\begin{lem} \label{lem04}
If $M$  is an $R$-complex,
then the induced morphisms 
\begin{align*}
\LL{a}{\deval{a}{M}}& \colon\LL{a}{\rhom_R(\Comp{R}{a},M)}\to\LL{a}{M} \\
\RG{a}{\deval{a}{M}}& \colon\RG{a}{\rhom_R(\Comp{R}{a},M)}\to\RG{a}{M} 
\end{align*}
are isomorphisms in $\D(R)$.
In particular, if $\LL{a}{M}\not\simeq 0$ or $\RG{a}{M}\not\simeq 0$,  then
$\rhom_R(\Comp{R}{a},M)\not\simeq 0$.
\end{lem}

\begin{proof}
For the first isomorphism,  
it suffices to check that the morphism
$$
\rhom_R(\RG{a}{R},\rhom_R(\Comp{R}{a},M))\xra{\rhom_R(\RG{a}{R},\deval{a}{M})}\rhom_R(\RG{a}{R},M)$$
is an isomorphism 
in $\D(R)$; see~\ref{props01}\eqref{item01d}.
In the following commutative diagram 
$$\xymatrix{
\rhom_R(\RG{a}{R},\rhom_R(\Comp{R}{a},M)) \ar[r]^-{(1)}_-{\simeq}\ar[d]_{\rhom_R(\RG{a}{R},\deval{a}{M})} 
    & \rhom_R(\RG{a}{R}\lotimes_R\Comp{R}{a},M) \ar[d]^{\rhom_R(\RG{a}{R}\lotimes_R\cpln{a}{R},M)} \\
\rhom_R(\RG{a}{R},M) & \rhom_R(\RG{a}{R}\lotimes_RR,M) \ar[l]_{\simeq}
}$$
(1) is adjunction and $\cpln{a}{R}\colon R\to\Comp{R}{a}$ is the natural inclusion.
Since $\RG{a}{R}\otimes_R\cpln{a}{R}$ is an
isomorphism  by~\ref{props01}\eqref{item01g}, the same is true of
$\rhom_R(\RG{a}{R}\lotimes_R\cpln{a}{R},M)$.
The diagram implies that $\rhom_R(\RG{a}{R},\deval{a}{M})$ is an isomorphism.

For the second isomorphism, use the equivalence~\ref{props01}\eqref{item01c}
to see that the vertical maps in the next commutative diagram are isomorphisms
$$\xymatrix{
\RG{a}{\rhom_R(\Comp{R}{a},M)} \ar[rr]^-{\RG{a}{\deval{a}{M}}}
\ar[d]_{\RG{a}{\nu^{}_{\rhom_R(\Comp{R}{a},M)}}}^{\simeq} 
    && \RG{a}{M} \ar[d]^{\RG{a}{\nu^{}_M}}_{\simeq} \\
\RG{a}{\LL{a}{\rhom_R(\Comp{R}{a},M)}} \ar[rr]^-{\RG{a}{\LL{a}{\deval{a}{M}}}}_-{\simeq} 
&& \RG{a}{\LL{a}{M}}.
}$$
The morphism $\RG{a}{\LL{a}{\deval{a}{M}}}$ is an isomorphism in $\D(R)$ because we have
shown that
$\LL{a}{\deval{a}{M}}$ is so.  The diagram shows that $\RG{a}{\deval{a}{M}}$ is an isomorphism
as well.

The final statement follows from
the additivity of $\LL{a}{-}$ and $\RG{a}{-}$.
\end{proof}

\begin{lem} \label{lem06}
If $M,N$ are homologically finite $R$-complexes, 
then the complex $X=\rhom_R(N,M)$ is homologically degreewise finite and 
$\LL{a}{X}\simeq X\lotimes_R{\Comp{R}{a}}$.  In particular, one has
$\inf(\LL{a}{X})=\inf(X)$ and 
$\sup(\LL{a}{X})=\sup(X)$.
\end{lem}

\begin{proof}
The finiteness of each $\HH_i(X)$ is standard.
A verification of the isomorphism 
is essentially in~\cite[Proof of (5.9)]{christensen:ogpifd}.
The flatness of $R\to\Comp{R}{a}$ implies
$\HH_i(\LL{a}{X})\cong\HH_i(X)\lotimes_R{\Comp{R}{a}}$,
and the equalities follow from the faithful flatness of $R\to\Comp{R}{a}$.
\end{proof}

We next prove a vanishing result akin to~\cite[(2.3)]{buchweitz:psrp}.  
Note that $M$ is not assumed to be finitely generated.

\begin{prop} \label{prop01}
Let $M$ be an 
$R$-module such that
the morphism $\nu^{}_M\colon M\to\LL{a}{M}$
is an isomorphism in $\D(R)$.
Then
$\ext^i_R(\Comp{R}{a},M)=0$ for each $i\neq 0$ and 
the evaluation map $\eval{a}{M}\colon\Hom_R(\Comp{R}{a},M)\to M$
is an isomorphism.
\end{prop}

\begin{proof}
Because the morphism $\nu^{}_M\colon M\to\LL{a}{M}$
is an isomorphism in $\D(R)$, the same is true of
$\rhom_R(X,\nu^{}_M)\colon\rhom_R(X,M)\to \rhom_R(X,\LL{a}{M})$
for each $R$-complex $X$.
From~\ref{props01}\eqref{item01e} it follows that
the morphism 
$$\rhom_R(\gamma^{}_X,M)\colon\rhom_R(X,M)\to \rhom_R(\RG{a}{X},M)$$
is an isomorphism in $\D(R)$.

The naturality of $\gamma$ provides the following commutative diagram in $\D(R)$
$$\xymatrix{
\RG{a}{R} \ar[rr]^-{\RG{a}{\nu^{}_R}}_-{\simeq} \ar[d]_{\gamma^{}_R} && \RG{a}{\LL{a}{R}} 
  \ar[d]^{\gamma^{}_{\LL{a}{R}}} \\
R \ar[rr]^-{\nu^{}_R} && \LL{a}{R}
}
$$
and an application of $\rhom_R(-,M)$ yields the next commutative diagram in $\D(R)$
$$\xymatrix{
\rhom_R(\LL{a}{R},M) \ar[rrr]^-{\rhom_R(\nu^{}_R,M)} 
\ar[d]_{\rhom_R(\gamma^{}_{\LL{a}{R})},M)}^{\simeq}
&&& \rhom_R(R,M) \ar[d]^{\rhom_R(\gamma^{}_R,M)}_{\simeq} \\
\rhom_R(\RG{a}{\LL{a}{R}} ,M) \ar[rrr]^-{\rhom_R(\RG{a}{\nu^{}_R},M)}_-{\simeq}
&&& \rhom_R(\RG{a}{R},M)
}
$$
where the vertical morphisms are isomorphisms because of the argument of the previous paragraph.
Hence, the morphism $\rhom_R(\nu^{}_R,M)$ is also an isomorphism.

Consider next the commutative triangle
$$\xymatrix{
R \ar[d]_{\nu^{}_R} \ar[rd]^{\cpln{a}{R}} \\
\LL{a}{R} \ar[r]^{\kappa}_{\simeq} & \Comp{R}{a}
}
$$
where $\kappa$ is gotten by taking degree 0 homology;
see~\ref{props01}\eqref{item01f}.  Apply 
$\rhom_R(-,M)$ to produce the next commutative diagram in $\D(R)$
$$\xymatrix{
\rhom_R(\Comp{R}{a},M) \ar[rrrd]^{\quad\rhom_R(\cpln{a}{R},M)} \ar[d]_{\rhom_R(\kappa,M)}^{\simeq} \\
\rhom_R(\LL{a}{R},M) \ar[rrr]^{\rhom_R(\nu^{}_R,M)}_{\simeq} &&& \rhom_R(R,M)
}
$$
which implies that $\rhom_R(\cpln{a}{R},M)$ is an isomorphism in $\D(R)$.

In the final commutative diagram
$$\xymatrix{
\rhom_R(R,M) \ar[rrd]^{\xi}_{\simeq} \ar[d]_{\rhom_R(\cpln{a}{R},M)}^{\simeq} \\
\rhom_R(\Comp{R}{a},M) \ar[rr]_-{\deval{a}{M}} && M
}
$$
the morphism $\xi$ is the natural evaluation isomorphism.  The diagram shows that
$\deval{a}{M}$ is an isomorphism in $\D(R)$.
Since $M$ is a module, this implies $\ext^i_R(\Comp{R}{a},M)=0$ for each $i\neq 0$
and further that the induced map 
$\HH_0(\deval{a}{M})\colon \Hom_R(\Comp{R}{a},M)\to M$ is bijective.
The definitions yield an equality $\HH_0(\deval{a}{M})=\eval{a}{M}$, 
completing the proof.
\end{proof}

\begin{disc} \label{prop03}
If $M$ is an $R$-module such that
$M\cong\Comp{M}{a}$, then $M\simeq\LL{a}{M}$. 
Indeed, the isomorphism $M\cong\Comp{M}{a}$ shows that $M$ is an $\Comp{R}{a}$-module.
Let $P$ be an $\Comp{R}{a}$-free resolution of $M$.  Then $P$ is an $R$-flat resolution
of $M$ consisting of $\fa$-adically complete modules.  Thus, one has
$\LL{a}{M}\simeq \Lambda^{\fa}(P)\cong P\simeq M$. 
\end{disc}

We are now in a position to give a useful alternate description of 
the analytic conductor submodule
$\fc{a}{M}$; see~\ref{maps} for the definitions of the maps.

\begin{prop} \label{CisHom}
Let $M$ be a finitely generated $R$-module. The homomorphisms
$\eff{a}{M}\colon\fc{a}{M}\to\Hom_R(\Comp{R}{a},M)$ and 
$\ceval{a}{M} \colon\Hom_R(\Comp{R}{a},M)\to\fc{a}{M}$ are inverse
isomorphisms.  In particular, $\Hom_R(\Comp{R}{a},M)$ is finitely generated over $R$.
\end{prop}

\begin{proof}
One checks from the definitions that the composition $\ceval{a}{M}\eff{a}{M}$
is the identity on $\fc{a}{M}$.  Hence, the first conclusion will be verified once
we show that $\ceval{a}{M}$ is bijective; the second conclusion will then follow,
as $\fc{a}{M}$ is finitely generated over $R$.

The module $\fc{a}{M}$ is $\fa$-adically complete, so Proposition~\ref{prop01}
implies that the evaluation map $\eval{a}{\fc{a}{M}}\colon\Hom_R(\Comp{R}{a},\fc{a}{M})\to \fc{a}{M}$
is bijective.
By Lemma~\ref{lem02} the map 
$\Hom_R(\Comp{R}{a},\inc{a}{M})\colon\Hom_R(\Comp{R}{a},\fc{a}{M})
\to\Hom_R(\Comp{R}{a},M)$ is an isomorphism.  
In particular, the composition
$\ceval{a}{M}=\eval{a}{\fc{a}{M}}\circ\Hom_R(\Comp{R}{a},\inc{a}{M})^{-1}$
is bijective, as desired.
\end{proof}

\section{Detecting completeness} \label{sec02}

\begin{para} \label{proofA}
\emph{Proof of  Theorem~\ref{thma}.}
The implication
\eqref{thmaitem1}$\implies$\eqref{thmaitem2} follows from
Proposition~\ref{prop01} and Remark~\ref{prop03}, and
\eqref{thmaitem2}$\implies$\eqref{thmaitem3} is trivial.

For the implication \eqref{thmaitem3}$\implies$\eqref{thmaitem1},
set $S=R/\ann_R(M)$.  A result of Gruson and 
Raynaud~\cite[Seconde Partie, Thm.~(3.2.6)]{raynaud:cpptpm}, and 
Jensen~\cite[Prop.~6]{jensen:vl} provides the following bound on the
projective dimension of $\Comp{S}{a}$ as an $S$-module:
\begin{equation} \label{bound} \tag{$\ast$}
\pd_S(\Comp{S}{a})\leq\dim(S)=\dim_R(M).
\end{equation}
Consider the following sequence of isomorphisms in $\D(R)$:
\begin{align*}
\rhom_R(\Comp{R}{a},M)
&\simeq\rhom_R(\Comp{R}{a},\rhom_S(S,M)) \\
&\simeq\rhom_S(\Comp{R}{a}\lotimes_RS,M) \\
&\simeq\rhom_S(\Comp{S}{a},M) .
\end{align*}
The first isomorphism follows from the fact that $M$ is naturally an $S$-module.
The second is adjunction, and the third is standard as $S$ is 
finitely generated over $R$.
Combining~\eqref{bound} with the displayed isomorphisms, the 
assumption 
$\ext^i_R(\Comp{R}{a},M)=0$ for all $i =1,\ldots,\dim_R(M)$
implies $\ext^i_R(\Comp{R}{a},M)=0$ for all $i\neq 0$.

It follows that the natural map
$\lambda\colon \Hom_R(\Comp{R}{a},M)\to \rhom_R(\Comp{R}{a},M)$
is an isomorphism in $\D(R)$.
Proposition~\ref{CisHom} implies that the composition
$\lambda\circ\eff{a}{M}\colon \fc{a}{M} \to
\rhom_R(\Comp{R}{a},M)$ is also an isomorphism in $\D(R)$.
Because $M$ is finitely generated, the natural morphism
$\mu\colon\LL{a}{M}\to\Comp{M}{a}$ is also an isomorphism in $\D(R)$.
These data yield the following commutative diagram
$$\xymatrix{\fc{a}{M}\ar[r]^-{\nu^{}_{\fc{a}{M}}}_-{\simeq} \ar[d]_{\inc{a}{M}} & \LL{a}{\fc{a}{M}} 
\ar[rr]^-{\LL{a}{\lambda\circ\eff{a}{M}}}_-{\simeq}
&& \LL{a}{\rhom_R(\Comp{R}{a},M)} \ar[d]^{\LL{a}{\deval{a}{M}}}_{\simeq} \\
M\ar[r]^-{\cpln{a}{M}} &\Comp{M}{a} && \LL{a}{M} \ar[ll]_-{\mu}^-{\simeq}
}$$
one sees that the composition of natural maps
$\fc{a}{M}\xra{\inc{a}{M}}M\xra{\cpln{a}{M}}\Comp{M}{a}$ is bijective.
Since $\cpln{a}{M}$ is also injective, the result now follows. \qed
\end{para}

\begin{disc} \label{referee}
As the referee indicated, 
one can interpret Theorem~\ref{thma}
as a statement about the
$\fa$-adic completeness of  $R/\ann_R(M)$ because
$M$ is $\fa$-adically complete
if and only if $R/\ann_R(M)$ is $\fa$-adically complete.
For the sake of completeness, we include a sketch of the proof.

For one implication, assume that $M$ is $\fa$-adically complete.
For each prime $\fp\in\ass_R(M)$, the injection $R/\fp\hookrightarrow M$
and the completeness of $M$ imply that $R/\fp$ 
is $\fa$-adically complete.
In particular, this is true for each minimal prime $\fp$ containing
$\ann_R(M)$, and it follows that the same is true for each
non-minimal prime $\fp$ containing
$\ann_R(M)$.  A prime filtration argument applied to $R/\ann_R(M)$
shows that $R/\ann_R(M)$ is $\fa$-adically complete.

Conversely, if $R/\ann_R(M)$ is $\fa$-adically complete, then
there exists an integer $r$ and a surjection $(R/\ann_R(M))^r\twoheadrightarrow M$,
and it follows that $M$ is $\fa$-adically complete.

From this fact, one  easily deduces the following: 
When $N$ is a second finitely generated $R$-module,
if $M$ is $\fa$-adically complete and $\supp_R(N)\subseteq\supp_R(M)$,
then $N$ is $\fa$-adically complete.
\end{disc}

\begin{para} \label{proofreferee2}
\emph{Proof of Theorem~\ref{referee2}.}
One implication is trivial.  For the other, assume that
$\Comp{R}{a}$ is module-finite over $R$.
As $\Comp{R}{a}$ is flat and module-finite over $R$, it is projective,
and so
$\ext^i_R(\Comp{R}{a},R)=0$ for each $i\neq 0$.
The completeness of $R$ follows from 
Theorem~\ref{thma}. \qed
\end{para}

The next example shows that the nontrivial implication in Corollary~\ref{referee2} 
fails if $\fa$ is not assumed to be in the Jacobson radical of $R$.

\begin{ex} \label{fail}
Let $k$ be a field and set $R=k\times k$ and $\fb=k\times 0$.  
The Jacobson radical of $R$ is $0$.
One checks readily
that $\Comp{R}{b}\cong 0\times k$, showing that
$R$ is not $\fb$-adically complete even though
$\Comp{R}{b}$ is module-finite over $R$.
\end{ex}

Theorem~\ref{thma}
provides the converse to~\cite[(2.3)]{buchweitz:psrp} when 
$R$ is local and $M$ is finitely generated.  
This is the implication \eqref{item02iii}$\implies$\eqref{item02i} in the next result.
The implication  \eqref{item02i}$\implies$\eqref{item02ii} is in~\cite[(3.7)]{frankild:volh}
or~\cite[(2.3)]{buchweitz:psrp}, while 
the implication \eqref{item02ii}$\implies$\eqref{item02iii} is trivial.

\begin{cor} \label{cor02}
Let $(R,\m)$ be a local ring.  For a finitely generated $R$-module $M$ 
the following conditions are equivalent.
\begin{enumerate}[\quad\rm(i)]
\item \label{item02i} $M$ is  $\m$-adically complete.
\item \label{item02ii} For each flat $R$-module $B$ and each $i\neq 0$, one has
$\ext^i_R(B,M)=0$.
\item \label{item02iii}
For each $i\neq 0$, one has
$\ext^i_R(\Comp{R}{m},M)=0$. \hfill\qed
\end{enumerate}
\end{cor}

With Theorem~\ref{thma} and Corollary~\ref{cor02} in mind, one may ask what the
finitely generated complete $R$-modules look like, say, when $R$ is not complete.  
Examples include the modules of finite length.  
We observe next that one can have complete $R$-modules of infinite length.

\begin{ex}
Let $(S,\n)$ be a non-Artinian complete local ring. Set
$R=S[X]_{(\n,X)}$ with maximal ideal
$\m=(\n,X)R$.  The ring $R$ is not $\m$-adically complete,
while the module $R/(X)R\cong S$ is $\m$-adically complete
and has infinite length.
\end{ex}

A finitely generated $R$-module $C$ is semidualizing if $R\xra{\simeq}\rhom_R(C,C)$.

\begin{cor} \label{cor01}
If $C$ is a semidualizing $R$-module
such that $\ext^i_R(\Comp{R}{a},C)=0$ for all $i\neq 0$,
then $R$ is $\fa$-adically complete.
\end{cor}

\begin{proof}
Theorem~\ref{thma} implies that $C$ is $\fa$-adically complete and hence 
$C\simeq C\otimes_R\Comp{R}{a}\simeq C\lotimes_R\Comp{R}{a}$.
By~\cite[(5.8)]{christensen:scatac} the complex $C\lotimes_R\Comp{R}{a}$
is $\Comp{R}{a}$-semidualizing.
This provides (1)  in the following sequence
while (4) and (5) are by hypothesis
\begin{align*}
\Comp{R}{a}
&\stackrel{(1)}{\simeq}\rhom_{\Comp{R}{a}}(C\lotimes_R\Comp{R}{a},C\lotimes_R\Comp{R}{a})\\
&\stackrel{(2)}{\simeq}\rhom_R(C,\rhom_{\Comp{R}{a}}(\Comp{R}{a},C\lotimes_R\Comp{R}{a})\\
&\stackrel{(3)}{\simeq}\rhom_R(C,C\lotimes_R\Comp{R}{a})\\
&\stackrel{(4)}{\simeq}\rhom_R(C,C)\\
&\stackrel{(5)}{\simeq} R
\end{align*}
(2) is adjunction~\cite[(1.5.2)]{christensen:scatac},
(3) is standard~\cite[(1.5.5)]{christensen:scatac}.
\end{proof}

Here is a version of Theorem~\ref{thma} for complexes.

\begin{prop} \label{prop05}
Let $M$ be a homologically degreewise  finite $R$-complex such that $\inf(\LL{a}{M})=\inf(M)$ and 
$\sup(\LL{a}{M})=\sup(M)$, e.g., if $M$ is homologically finite.
Fix an integer  $s\geq\sup(M)$.
If $\rhom_R(\Comp{R}{a},M)$ is homologically concentrated in degree $s$,
then so is $M$, and the module $\HH_s(M)$ is $\fa$-adically complete.
\end{prop}

\begin{proof}
Assume $M\not\simeq 0$.  Then $\sup(\LL{a}{M})=\sup(M)>-\infty$, and 
Lemma~\ref{lem04} implies $\rhom_R(\Comp{R}{a},M)\not\simeq 0$.
Our hypotheses provide (1) and (3) in the sequence
$$
s
\stackrel{(1)}{\geq}\sup(M)
\stackrel{(2)}{\geq}\sup(\rhom_R(\Comp{R}{a},M))
\stackrel{(3)}{=}s
$$ 
and (2) is from~\cite[(2.1)]{foxby:ibcahtm};
this implies $s=\sup(M)$.  
Since
$\rhom_R(\Comp{R}{a},M)$ is homologically concentrated in degree $s$,
one has 
$\shift^s\ext^{-s}_R(\Comp{R}{a},M)\simeq\rhom_R(\Comp{R}{a},M)$, providing the first
of the following isomorphisms
$$\LL{a}{\shift^s\ext^{-s}_R(\Comp{R}{a},M)}
\simeq\LL{a}{\rhom_R(\Comp{R}{a},M)}
\simeq\LL{a}{M}$$
while the second one is from Lemma~\ref{lem04}.
This provides (5) in the next sequence
$$
\inf(M)
\stackrel{(4)}{=}\inf(\LL{a}{M})
\stackrel{(5)}{=}\inf(\LL{a}{\shift^s\ext^{-s}_R(\Comp{R}{a},M)})
\stackrel{(6)}{\geq} s
\stackrel{(7)}{=}\sup(M)
\stackrel{(8)}{\geq}\inf(M)
$$
while (4) is by assumtion,
(6) is by~\ref{props01}\eqref{item01j},
(7) is proved above,
and (8) is trivial.
It follows that $\inf(M)=\sup(M)=s$ and so $M$ is homologically concentrated in
degree $s$.  Finally, one has $M\simeq\shift^s\HH_s(M)$ and so 
$$\rhom_R(\Comp{R}{a},M)
\simeq\rhom_R(\Comp{R}{a},\shift^s\HH_s(M))
\simeq\shift^s\rhom_R(\Comp{R}{a},\HH_s(M)).$$
Since this is homologically concentrated in degree $s$, one has
$\ext^R_i(\Comp{R}{a},\HH_s(M))=0$ for each $i\neq 0$.  
Theorem~\ref{thma} implies that $\HH_s(M)$ is $\fa$-adically complete.
\end{proof}

The next result contains Theorem~\ref{thmb} from the introduction.

\begin{cor} \label{cor03}
Let $M,N$ be a homologically finite $R$-complexes
and $s\in\mathbb{Z}$  such that $s\geq\sup(\rhom_R(N,M))$.
If $\rhom_R(\Comp{N}{a},M)$ is homologically concentrated in degree $s$,
then so is $\rhom_R(N,M)$, and $\ext_R^{-s}(N,M)$ is $\fa$-adically complete.
\end{cor}

\begin{proof}
\ref{props01}\eqref{item01f} and adjunction provide the following sequence.
$$\rhom_R(\Comp{N}{a},M)\simeq\rhom_R(\Comp{R}{a}\lotimes_RN,M)\simeq
\rhom_R(\Comp{R}{a},\rhom_R(N,M))$$
Lemma~\ref{lem06} shows that Proposition~\ref{prop05} applies to the complex
$\rhom_R(N,M)$, yielding the desired conclusion.  
\end{proof}

\begin{cor} \label{cor04}
Assume that $R$ is local and $M,N$ are nonzero 
finitely generated $R$-modules with $\pd_R(N)<\infty$.  
If $\ext_R^i(\Comp{N}{a},M)=0$ for each $i\neq 0$, then $N$ is free and $M$ is $\fa$-adically 
complete.
\end{cor}

\begin{proof}
Using $s=0$ in Corollary~\ref{cor03}, one concludes that
$\ext_R^i(N,M)=0$ for each $i\neq 0$ and that $\Hom_R(N,M)$ is $\fa$-adically complete.
Since $\pd_R(N)$ is finite, one has
$$\rhom_R(N,M)\simeq\rhom_R(N,R)\lotimes_R M$$
by tensor-evaluation~\cite[(4.4)]{avramov:hdouc}.  
The next equalities are from~\cite[(2.1)]{foxby:ibcahtm}
and~\cite[(2.13)]{christensen:scatac}.
$$0=\inf(\rhom_R(N,M))=\inf(\rhom_R(N,R))+\inf(M)=-\pd_R(N)$$
Since $R$ is local, the module $N\neq 0$ is free and 
$\Hom_R(N,M)\cong M^n$ for some $n>0$.
Because $M^n$ is
$\fa$-adically complete, the same is true of $M$.
\end{proof}

\section*{Acknowledgments}

We are grateful to Phillip Griffith, Srikanth Iyengar, 
Christian U.~Jensen, and Anders Thorup for stimulating discussions about 
this research.  We also thank the anonymous referee for helpful comments.

\providecommand{\bysame}{\leavevmode\hbox to3em{\hrulefill}\thinspace}
\providecommand{\MR}{\relax\ifhmode\unskip\space\fi MR }
\providecommand{\MRhref}[2]{%
  \href{http://www.ams.org/mathscinet-getitem?mr=#1}{#2}
}
\providecommand{\href}[2]{#2}

\end{document}